\documentclass{elsart3-1}

\usepackage{amssymb}
\usepackage{amsmath}
\usepackage[english,francais]{babel}

\newtheorem{theoreme}{Th\'eor\`eme}[section]
\newtheorem{lemme}[theoreme]{Lemme}
\newtheorem{proposition}[theoreme]{Proposition}
\newtheorem{corollaire}[theoreme]{Corollaire}

{\nopagebreak[0]\hfill$\Box$ \medskip}%
\renewcommand{\theequation}{\arabic{equation}}
\def\Q{{\bf Q}}

\def\F{{\bf F}}
\def\Z{{\bf Z}}

\def\t{{\theta}}

\def\m{{\mathfrak m}}
\def\n{\noindent}
\def\dim{\mathop{\rm dim}}
\def\End{\mathop{\rm End}}
\def\ar{\Longleftrightarrow}
\def\cf{\mathcal{F}}

\def\D{\Delta}

\usepackage[T1]{fontenc}

\def\og{\leavevmode\raise.3ex\hbox{$\scriptscriptstyle\langle\!\langle$~}}
\def\fg{\leavevmode\raise.3ex\hbox{~$\!\scriptscriptstyle\,\rangle\!\rangle$}}

\journal{the Acad\'emie des sciences}
\begin{document}

\centerline{}
\begin{frontmatter}




%
\selectlanguage{francais}
\title{Formes modulaires modulo $2:$ structure de l'algèbre de Hecke}



\author[authorlabel1]{Jean-Louis NICOLAS},
\ead{jlnicola@in2p3.fr,\;\;http://math.univ-lyon1.fr/$\sim$nicolas/}
\author[authorlabel2]{Jean-Pierre SERRE}
\ead{jpserre691@gmail.com}

\address[authorlabel1]{CNRS, Universit\'e de Lyon, Institut Camille
  Jordan, Math\'ematiques,  F-69622 Villeurbanne Cedex, France.}
\address[authorlabel2]{Coll\`ege de France, 3 rue d'Ulm, F-75231 Paris Cedex 05, France.}


\medskip

\selectlanguage{english}
\begin{abstract}



\noindent
{\bf Modular forms mod $2 :$ structure of the Hecke ring}

\vskip 0.5\baselineskip

We show that the Hecke algebra for modular forms mod 2 of level~1 is isomorphic to the power series ring $\F_2[[x,y]]$, where $x=T_3$
and $y = T_5$.

\medskip

\n

\medskip

\n
\end{abstract}
\end{frontmatter}


\selectlanguage{francais}

\section{Notations}
\label{par1}

  Nous conservons les notations de la Note précédente \cite{NS}. En particulier,
  nous notons $\D $ l'élément de $\F_2[[q]]$ défini par:
  
  \smallskip
  
  \quad $ \D= \sum_{n=1}^\infty \tau(n)q^n =  \sum_{m=1}^\infty q^{(2m+1)^2},$ 
  
  \smallskip
  
  \n  et \  $\cf$ désigne le sous-espace vectoriel de $\F_2[[q]]$ engendré par les puissances impaires de $\D$  :
    \smallskip
  
 \quad  $\cf = \ < \! \! \D, \D^3, \D^5, ... \! \!>.$
  
  
  
    
   
   \smallskip
   
   \n L'espace
   $\cf$ est stable par les opérateurs de Hecke $T_p$ , 
   $p$ premier $\neq 2.$

\section{Les espaces $\cf(n)$ et les algèbres $A(n)$}
\label{par2}
   
Soit $n$ un entier $> 0$. Soit $\cf(n)$ le sous-espace de $\cf$ de base
$\{\D, \D^3, ..., \D^{2n-1}\}$. On a  $\dim \cf(n) = n$.
Soit $A(n)$ la $\F_2$-sous-algèbre de $\End(\cf(n))$ engendrée par $\F_2$ et les $T_p$. On a $A(n) = \F_2 \oplus \m (n)$, où  $\m(n)$  est l'unique idéal maximal de $A(n)$ (à savoir le sous-espace vectoriel de $A(n)$ engendré par les $T_p$ et leurs produits); cet idéal est nilpotent.

  Soit $\cf(n)^*$ le dual de l'espace vectoriel $\cf(n)$, muni de sa structure 
naturelle de  $A(n)$-module, et soit $e_n$ l'élément de $\cf(n)^*$ défini par:
  
  \smallskip 
  
 \quad     $< \! e_n,\D \!> \ = 1$  \quad  et  \quad $ <\! \!e_n,\D^{2i+1} \! \! >  \ = 0$ \ si  \ $1 \leqslant i < n$.
     
      \smallskip
      
     \n Si $f = \sum a_m(f) q^m$  est un élément de $\cf(n)$, on a:
     
      \smallskip
     
 \quad     $<\! \!e_n,f\! \!> \ = \ a_1(f)$ \quad et \quad  $ < \! T_pe_n,f\! \!> \ = \ a_p(f)$ \  pour tout  $p$. 
 
  \smallskip
     
   \n  On en déduit par récurrence sur  $r$  la formule:
   
\smallskip

\begin{equation}\label{2.1}     
<  \!T_{p_1}...T_{p_r}e_n,f \!> \ = \ a_{p_1...p_r}(f),
\end{equation}
 
\smallskip

\n où les  $p_i$ sont des nombres premiers $\neq2$.

\smallskip

\begin{lemme}\label{lem1}
{\bf -}
Soit $f \in \cf(n), f \neq 0$. Il existe $u \in A(n)$ 
tel que \ $<\! \!e_n,uf\! \!> \ = \ 1$.
\end{lemme}

\smallskip

\n {\it Démonstration}. Ecrivons $f$ sous la forme $f = q^m + \sum_{i>m}a_iq^i$ et soit $m = p_1...p_r$ une décomposition de $m$ en produit de nombres premiers. Comme  $m$  est impair, il en est de même des $p_i$. Soit $u = T_{p_1}...T_{p_r}.$ La formule (\ref{2.1})  montre que
$<\! \!ue_n,f\! \!> = 1$. Comme $<\! \!ue_n,f\! \!> = <\! \!e_n,uf\! \!>$, cela démontre le lemme.

\section{Quelques propriétés de $\cf(n)$ et de $A(n)$}
\label{par3}

\begin{proposition}\label{prop1}
{\bf -}
Le $A(n)$-module $\cf(n)^*$ est libre de base \ $e_n$.
\end{proposition}
\smallskip

\n {\it Démonstration}. Soit $E$ le sous-$A(n)$-module de $\cf(n)^*$ engendré par  $e_n$. Si $E$ était distinct de $\cf(n)^*$, il existerait $f \in \cf(n), f \neq 0$, tel que \ $<\! \!ue_n,f\! \!> = 0$ pour tout $u \in A(n)$, ce qui 
contredirait le lemme 1. On a donc $E = \cf(n)^*$, ce qui montre que
$\cf(n)^*$ est engendré par $e_n$. D'où la proposition.


\smallskip

\n {\it Remarque}. Si $n > 2$, le $A(n)$-module $\cf(n)$ n'est pas un module libre.

\smallskip

\smallskip

\begin{corollaire}\label{coro1}
{\bf -}
L'application $A(n) \to \cf(n)^*$ donnée par $u \mapsto ue_n$ est bijective.
\end{corollaire}

\smallskip

  Ce n'est qu'une reformulation de la proposition. Noter que, par dualité, on obtient ainsi une bijection de $\cf(n)$ sur le dual $A(n)^*$ de l'espace vectoriel $A(n)$.
  
\smallskip

\begin{corollaire}\label{coro2}
{\bf -}
On a $\dim A(n) = n$.
\end{corollaire}

  \smallskip
  
  Cela résulte du Corollaire pr\'ec\'edent et du fait que $\dim \cf(n) = n$.
  
  \smallskip

\begin{corollaire}\label{coro3} 
{\bf -}
Le commutant de $A(n)$ dans $\End(\cf(n))$ est égal à $A(n)$.
\end{corollaire}

\smallskip

  Par dualité, cela revient à dire que le commutant de $A(n)$ dans $\End(\cf(n)^*)$ est égal à $A(n)$, ce qui résulte de la proposition.
  
\smallskip
  
\begin{proposition}\label{prop2}
{\bf -}
L'algèbre $A(n)$ est engendrée par $T_3$ et $T_5$.
\end{proposition}

\smallskip

\n {\it Démonstration}. Soit $A'$ la sous-algèbre de $A(n)$ engendrée par $T_3$ et $T_5$. C'est une algèbre locale; soit $\m'$ son idéal maximal. Supposons que $A' \neq A(n)$, i.e. $\dim A' < n$. Le $A'$-module $\cf(n)^*$
n'est pas monogène: sinon, sa dimension serait $<n$. D'après le lemme de Nakayama, cela signifie que le quotient $V = \cf(n)^*/\m'\cf(n)^*$ est de dimension $> 1$. Par dualité, cela équivaut à dire que le sous-espace de
$\cf(n)$ annulé par $\m'$ est de dimension $>1$. Il existe donc $f \in \cf(n)$,
avec $f \neq 0, \D$, tel que $T_3f = T_5f = 0$, ce qui contredit le
cor.5.3 au th.5.1 de \cite{NS}.

\section{Passage à la limite: l'algèbre $A$}
\label{par4}

On a $\cf(n) \subset \cf(n+1)$ et la restriction à $\cf(n)$ d'un élément de $A(n+1)$ appartient à $A(n)$. On obtient ainsi un homomorphisme surjectif
$A(n+1) \to A(n)$. D'où un système projectif

  $ ... \to A(n+1) \to A(n) \to \ \ ... \ \ \to A(2) \to A(1) =\F_2.$ 

Nous noterons $A$ la limite projective de ce système. C'est un anneau
local commutatif; il est compact pour la topologie limite
projective. Son idéal maximal $\m$ est la limite projective des $\m(n)$.  
L'anneau $A$ opère de façon naturelle sur $\cf$.

Soient  $x$  et $y$  deux indéterminées. Pour chaque  $n$, il existe 
un unique homomorphisme $\psi_n:\F_2[x,y] \to A(n)$  tel que
$\psi_n(x) = T_3$ et $\psi_n(y) = T_5$. Par passage à la limite, 
on en déduit un homomorphisme 

  $\psi : \F_2[[x,y]] \ \to \ A$  

\n tel que $\psi(x) = T_3$ et $\psi(y) = T_5$.

\smallskip

\begin{theoreme}\label{thm1}
{\bf -}
L'homomorphisme $\psi$ défini ci-dessus est
un isomorphisme.
\end{theoreme}

\n {\it Démonstration}. La surjectivité de $\psi$ résulte de la prop. \ref{prop2}. Pour prouver l'injectivité, il suffit de montrer que, pour tout élément  $u = \sum \lambda_{ij}x^iy^j$ non nul de $\F_2[[x,y]]$, il existe $f\in \cf$ tel que:

\smallskip

\begin{equation}\label{4.1}
\sum \lambda_{ij}T_3^i T_5^j f = \D.
\end{equation}

\smallskip

\n [Noter que la somme est une somme finie, car $T_3^i T_5^j f = 0$  quand  $i+j$  est assez grand (par exemple $i+j > \deg f$).]

\smallskip

Si $\lambda_{00} = 1$ on prend  $f=\D$. Supposons donc $\lambda_{00} = 0$. Soit $\Sigma$ l'ensemble des couples $(i,j)$ avec $\lambda_{ij} =1$; considérons ceux pour lesquels l'entier $i+j$ est minimal, et parmi ceux-là, soit $(a,b)$ le couple où $a$ est maximum.  Soit  $k$  l'entier impair de code $[a,b]$, au sens de \cite[\S4.1]{NS} et soit $f = \D^k$. On montre, en utilisant les Propositions 4.3 et 4.4 de \cite[\S4]{NS}, que l'on a
$T_3^aT_5^bf = \D$ et $T_3^iT_5^jf = 0$ pour tout $(i,j) \in \Sigma$ distinct de $(a,b)$. D'où (\ref{4.1}).

\smallskip

\begin{corollaire}\label{coro4}
{\bf -}
L'algèbre $A$ est un anneau local régulier de
dimension $2$. En particulier, c'est un anneau intègre.
\end{corollaire}

\smallskip

  A partir de maintenant, nous identifierons les algèbres $A$ et $\F_2[[x,y]]$
  au moyen de $\psi$; cela nous permettra d'écrire  $x$  et  $y$  à la place de $T_3$ et $T_5$.

\section{Structure des  $A$-modules $\cf$ et $\cf^*$}
\label{par5}

\n    L'algèbre  $A$  opère sur $\cf$. Par dualité, elle opère aussi sur le dual
    $\cf^*$ de $\cf$, qui  est la limite projective des $\cf(n)^*$. Soit $e \in \cf^*$ la forme linéaire sur $\cf$ définie par:
      
      \smallskip
      
     \n  $< \! e,f \!> = a_1(f)$  pour tout  $f\in \cf$,  où  $a_1(f)$ désigne le coefficient de $q$ dans $f$.

\smallskip
      
\begin{theoreme}\label{thm2}
{\bf -}
 a)  Le $A$-module $\cf^*$ est libre de base $e$.
      
b) Le $A$-module $\cf$ est isomorphe à l'espace $A^*_{\rm cont}$ des
formes linéaires continues sur $A$.
\end{theoreme}  

\smallskip
   
\n[Une forme linéaire sur $A$ est continue si et seulement si
      elle s'annule sur une puissance de l'idéal maximal de $A$.]

 \n {\it Démonstration}. L'assertion a) résulte de la prop.\ref{prop1} par dualité; il en est de même de b) car $A^*_{\rm cont} = \cup_{n  \geqslant 1}  A(n)^*$.

\smallskip

\begin{corollaire}\label{coro5}
{\bf -}
Le $A$-module $\cf$ est divisible$:$ pour tout $u\in A,
u \neq 0$, la multiplication par $u$ est un endomorphisme surjectif de $\cf$.
En particulier, les endomorphismes $T_p: \cf \to \cf$ sont surjectifs.
\end{corollaire}

\smallskip

\n {\it Démonstration}. Par dualité, cela revient à dire que $u : \cf^* \to \cf^*$
est injectif, ce qui est clair puisque $A$ est un anneau intègre.

\n{\it Remarque.} D'après \cite{No}, $\cf$ est un $A$-module {\it injectif}, à savoir
l'enveloppe injective du corps résiduel $\F_2$ de $A$. C'est là une propriété plus forte que la propriété de divisibilité.

\section{Une base de $\cf$ adaptée à $T_3$ et $T_5$}
\label{par6}

\begin{theoreme}\label{thm3}
{\bf -}
Il existe une base $m(a,b)_{a,b \geqslant 0}$ de $\cf$ et une seule 
qui a les quatre propriétés suivantes$~ :$

i) $m(0,0) = \D.$

ii) $ < \! e,m(a,b) \! > \ = \ 0 $ \ {\it si} \ $a+b > 0$.

iii) $ T_3|m(a,b) = \left\{
\begin{array}{ll}
\! \! m(a \! -1,b) \ \ \ {\it si}  \ a > 0 \\

\! \!  0  \quad \quad \quad \quad \quad {\it si}  \  a = 0.

\end{array}
\right.
$

iv) $ T_5|m(a,b) = \left\{
\begin{array}{ll}
 \! \! m(a,b \! - \! 1) \ \ \ {\it si}   \ b > 0 \\

\! \! 0  \quad \quad \quad \quad \quad {\it si} \  b = 0.

\end{array}
\right.
$
\end{theoreme}

\n {\it Démonstration}. D'après le th.\ref{thm2}, il suffit de prouver le même énoncé pour le $A$-module $A^*_{\rm cont}$, et dans ce cas on définit $m(a,b)$ comme la forme linéaaire sur $A$ donnée par:

\smallskip

$  \quad \sum n_{ij}x^iy^j \ \mapsto  \ n_{ab}.$

\smallskip

\n Les propriétés i) à iv) sont évidentes. L'unicité se démontre par récurrence sur $a+b$.

\smallskip

\n {\it Exemples} (cf. \cite{siteweb}):

\n  $ m(0,0)= \D; \ m(1,0)= \D^3; \ m(0,1)=\D^5;$ 

\smallskip

\n $m(2,0)= \D^9; \ m(1,1)=\D^7; \ m(0,2)=\D^{17}; $

\smallskip 

\n $m(3,0)=\D^{11}; \ m(2,1)=\D^{13};  \ m(1,2)= \D^{11}+\D^{19};  \ m(0,3)= \D^{13}+\D^{21};$

\smallskip

\n $m(2^r,0)=\D^{1+2^{2r+1}},\ m(2^r \! - \! 1,0)=  \D^{(1+2^{2r+1})/3}$ \ et \ \ $m(0,2^r)= \D^{1+2^{2r+2}}$.

\smallskip

\n {\it Remarques}. 

1)  {\it L'exposant dominant} de $m(a,b)$ au sens de \cite[\S4.3]{NS} est l'entier impair de code $(a,b)$; cela se déduit des résultats énoncés dans \cite[\S4]{NS}.
En particulier, l'ordre de nilpotence de  $m(a,b)$ est égal à $a+b+1$.

2) D'après Macaulay (\cite{Ma}, voir aussi \cite{No}) il est commode de noter les
$m(a,b)$ comme des monômes  $x^{-a}y^{-b}$, avec la convention que
$x^{-a}y^{-b}=0$ si $a$ ou $b$ est $<0$. Les formules du th.\ref{thm3} s'écrivent
alors simplement

\smallskip

  $ x.x^{-a}y^{-b} = x^{1-a}y^{-b}$ \ et \ $y.x^{-a}y^{-b} = x^{-a}y^{1-b}$.

\section{Développement des $T_p$ comme séries en $x=T_3$ et $y=T_5$}
\label{par7}

D'après le th.\ref{thm1}, tout $T_p$ peut s'écrire comme une série formelle en $x=T_3$ et $y=T_5$:
\begin{equation}\label{7.1}
T_p = \sum_{i+j\geqslant 1} a_{ij}(p)x^iy^j, {\rm \; \;avec\;\;} \ a_{ij}(p) \in \F_2.
\end{equation}

\smallskip

\n De façon plus précise, on a:
\begin{equation}\label{7.2}
T_p \in \F_2[[x^2,y^2]]  \ \ \ \ \;\;{\rm si\;}  p\equiv 1\pmod{8},
\end{equation}
\begin{equation}\label{7.3}
T_p \in x.\F_2[[x^2,y^2]]  \  \ \;\;  {\rm si\;}   p\equiv 3\pmod{8},
\end{equation}
\begin{equation}\label{7.4}
T_p \in y.\F_2[[x^2,y^2]]  \ \  \;\;{\rm si\;}  p\equiv 5\pmod{8},
\end{equation}
\begin{equation}\label{7.5}
T_p \in xy.\F_2[[x^2,y^2]]   \  \;\;{\rm si\;}  p\equiv 7\pmod{8}.
\end{equation}

\smallskip

\n {\it Exemples} (cf. \cite{siteweb}) :

$T_{17}=x^2+y^2+x^2y^2+x^6+x^4y^2+y^6+x^6y^2+x^4y^4+x^2y^6+
x^{10}+x^{10}y^2+x^6y^6+x^4y^8+x^2y^{10}+\ldots$

$T_{11}= x(1+x^2+y^2+x^4+x^2y^2+y^4+x^2y^4+y^6+x^6y^2+
x^8y^2+x^6y^4+x^2y^8+y^{10}+x^{10}y^2+\ldots)$

$T_{13}= y(1+x^2+y^2+x^4+y^4+x^6+x^4y^2+x^2y^4+x^6y^2+x^2y^6+
y^8+x^{10}+x^8y^2+x^6y^4+y^{10}+\ldots)$

 $T_7 = xy(1 + x^2 + x^4+x^2y^2+y^6+
x^6y^2+y^8+x^{10}+x^8y^2+x^6y^4+x^{12}+x^4y^8+x^2y^{10}+\ldots)$

\smallskip

\n Dans des cas simples, on peut donner explicitement la valeur du coefficient $ a_{ij}(p)$. Par exemple:
\begin{equation}\label{7.6}  
a_{10}(p) = 1  \ar \  p\equiv 3\;\pmod{8}
\end{equation}
\begin{equation}\label{7.7}  
a_{01}(p) = 1  \ar \  p\equiv 5\;\pmod{8}
\end{equation}
\begin{equation}\label{7.8}  
a_{11}(p) = 1  \ar \  p\equiv 7\pmod{16}
\end{equation}
\begin{equation}\label{7.9}  
a_{20}(p)= 1 \ar p {\rm \ est\; de\; la\; forme\;\;} a^2+8b^2 {\rm \;
  avec\;\;} a,b \in \Z, \; b {\rm \; impair}
\end{equation}
\begin{equation}\label{7.10} 
 a_{02}(p)= 1 \ar p {\rm\ est\; de\; la\; forme\;\;}  a^2+16b^2 {\rm
   \;avec\;\;} a,b \in \Z,  b {\rm \;impair}.
 \end{equation}
  
  \smallskip
  
 \n  Les formules (\ref{7.3}) et (\ref{7.6}) montrent que, si 
$ p\equiv 3\pmod{8}$, alors $T_p$  est le produit de $x$ 
par une série inversible en $x^2$ et $y^2$; en particulier, 
$T_p$ et $T_3$ ont le même noyau. Même chose si 
$ p\equiv 5\pmod{8}$ avec $x$ et $T_3$ remplacé par $y$ et $T_5$. On 
en déduit que l'algèbre $A$ est topologiquement engendrée par n'importe 
quel couple $(T_p, T_{p'}) $ avec $ p\equiv 3\pmod{8}$  et 
$ p'\equiv 5\pmod{8}$. Notons aussi que la proposition 4.3 (resp. 4.4) 
de \cite{NS} reste valable si l'on remplace 
$T_3$ par $T_p$ avec $p\equiv 3 \pmod{8}$ (resp. 
$T_5$ par $T_{p'}$ avec $p'\equiv 5 \pmod{8}$). 
  
  \smallskip
  
  \n {\it Remarques.} 

\smallskip

1) Pour  $i$ et $j$ fixés, la fonction $p \mapsto a_{ij}(p)$ 
est {\it frobénienne} au sens de
  \cite[\S3.3]{Se}. De façon plus précise, sa valeur ne dépend 
que de la substitution de Frobenius de $p$ dans une certaine 
extension galoisienne finie de  $\Q$, qui est non ramifiée 
en dehors de $\{2\}$ et dont le groupe de Galois est un $2$-groupe. 
Dans les deux premiers exemples ci-dessus, on peut prendre pour 
extension galoisienne le corps $\Q(\mu_8)$ des racines huitièmes 
de l'unité; dans les trois autres, les corps $\Q(\mu_8,\sqrt{uv})$, 
$\Q(\mu_8,\sqrt{u})$ et $\Q(\mu_8,\sqrt{v})$ avec $u=1+i$ et $v=\sqrt{2}$;
le premier de ces corps est le corps $\Q(\mu_{16})$ des racines $16$-i\`emes
de l'unit\'e; les deux autres ont des groupes 
de Galois sur $\Q$ qui sont diédraux d'ordre 8.

\smallskip

 2) Si $p > 5$, on peut se demander si la série donnant  $T_p$
peut \^etre un polyn\^ome en $x$ et $y$. La réponse est ``non'' :
d'apr\`es un r\'esultat r\'ecent de J. Bella\"iche, les $T_p$ sont
{\it algébriquement indépendants} sur $\F_2$.

\section{Séries thêta associées à $\Q(\sqrt{-2})$}
\label{par8}

  Soient $n$ un entier $\geqslant 1$ et soit $t \in \Z/2^n\Z$. Soit $\t_{t,n} \in \F_2[[q]]$ la série définie par:
  
  \smallskip
  
  \n $$  \t_{t,n} = \sum_{a \ {\rm impair} \, > \,  0} \ \ \sum_{b \equiv ta  \!\!\!\pmod{2^n}}q^{a^2+2b^2}.$$
  
  \smallskip
  
\n  On a:
  
  $\t_{0,n} = \D, \ \ \t_{t,n} = \t_{-t,n},  \ \ \t_{2^{n-1},n}= 0, \ \ \t_{t,n}+\t_{2^{n-1}-t,n} = \t_{t,n-1},$ 
  et \ \ $ \t_{2^{n-2},n} = \D^{1+2^{2n-3}}$ \ si $n \geqslant 2$.
  
  \smallskip
  
\n  Les séries $\t_{t,n}$ appartiennent à $\cf$. De façon plus précise:
 
\smallskip
 
\begin{theoreme}\label{thm4}
{\bf -}
Pour $n > 0$ fixé, les $\t_{t,n}$ engendrent le même sous-espace 
vectoriel de $\cf$ que les formes  $m(a,0)$ 
avec  $0 \leqslant a < 2^{n-1}$.
\end{theoreme} 

\smallskip
 
\n [Pour la définition des $m(a,b)$, voir \S\ref{par6}.]

\smallskip

\begin{corollaire}\label{coro6}
{\bf -}
Soit $f = \sum a_nq^n $ un élément de $\cf$. Les propriétés suivantes sont équivalentes  $:$

  $1) \ T_5|f = 0$.
  
  $ 2)$  La série $f$ est de la forme $\sum \t_{t_i,n_i}$.
   
   $3) \ a_n = 1 \ \ \Rightarrow \ \  n$ est de la forme $a^2+2b^2$, avec $a,b \in \Z$.
\end{corollaire}

\smallskip
  
\n {\it Exemples }(la table des $\t_{t,n}$ pour $n\leq 6$ et $0\leq t
\leq 2^{n-1}$ est sur le site \cite{siteweb}) :
  
 $\t_{0,1} = \D;$
 
 $\t_{0,2} = \D; \ \t_{1,2} = \D^3;$
 
 $\t_{0,3} = \D; \ \t_{1,3} = \D^3 + \D^{11}; \ \t_{2,3} = \D^9; \ \t_{3,3} = \D^{11}$.
 
 \medskip
 \n {\it Action des opérateurs de Hecke sur les $\t_{t,n}$}.
 
 Si $p \equiv 5 \ {\rm ou} \ 7 \pmod 8$, on a $T_p|\t_{t,n} = 0.$ 
 
 Si $p \equiv 1 \ {\rm ou} \ 3 \pmod 8$, on écrit  $p$  sous la forme
 $p = a^2 + 2b^2$, avec  $a,b \in \Z$; on définit $t(p) \in \Z/2^n\Z$ par
 $t(p) \equiv b/a  \pmod {2^n}$, et l'on pose $t^*(p) = -t(p)$.
   On a:
 
 \smallskip
 
  $T_p|\t_{t,n} = \t_{t\bullet t(p),n} + \t_{t\bullet t^*(p),n}$
 
 \smallskip
 
 \n où l'on a noté $x\bullet y$ la loi de composition\footnote{Cette 
loi munit  $\Z/2^n\Z$  d'une structure de groupe
abélien; ce groupe est cyclique d'ordre $2^n$; on peut l'interpréter 
comme le groupe des classes de formes quadratiques binaires 
primitives de discriminant $-2^{2n+3}$, ou encore comme le 
groupe Pic du sous-anneau de $\Z[\sqrt{-2}]$ de conducteur $2^n$.}
sur $\Z/2^n\Z$ définie par la formule $x\bullet y = (x+y)/(1-2xy)$.

\n On a en particulier $\;\t_{2^{n-1}-t(p),n} = T_p|\D^{1+2^{2n-1}}$. 
 
\section{Séries thêta associées à $\Q(i)$}
\label{par9}

  Les définitions et les résultats sont essentiellement les mêmes que ceux
  du \S\ref{par8}, à cela près que $a^2+2b^2, \ T_5$ et $m(a,0)$ sont remplacés par
   $a^2+4b^2, \ T_3$ et $m(0,b)$. De façon plus précise, si  $t$ et $n$ sont comme ci-dessus,  on définit la série thêta d'indice $(t,n)$ par:

\smallskip

\n $$  \t_{t,n}' = \sum_{a \ {\rm impair} \,  >  \, 0} \ \ \sum_{b \equiv ta \! \! \!\pmod{2^n}}q^{a^2+4b^2}.$$
  
\smallskip

\n On a:

\quad  $\t_{0,n}' = \D, \quad \t_{t,n}' = \t_{-t,n}', \quad
\t_{2^{n-1},n}'= 0, \quad \t_{t,n}'+\t_{2^{n-1}-t,n}' = \t_{t,n-1}', $
et \ \ $ \t_{2^{n-2},n}' = \D^{1+2^{2n-2}}$ \ si $n \geqslant 2 $.

\smallskip

\n De plus:

\smallskip

\begin{theoreme}\label{thm5}
{\bf -}
Pour $n > 0$ fixé, les $\t_{t,n}'$ engendrent le même sous-espace
vectoriel de $\cf$ que les formes  $m(0,b)$ avec  $0 \leqslant b < 2^{n-1}$.
\end{theoreme}

\smallskip

\begin{corollaire}\label{coro7}
{\bf -}
Soit $f = \sum a_nq^n $ un élément de $\cf$. Les propriétés suivantes sont équivalentes $:$

  $1) \ T_3|f = 0$.
  
  $ 2)$  La série $f$ est de la forme $\sum \t_{t_i,n_i}'$.
   
   $3) \ a_n = 1 \ \ \Rightarrow \ \ n$ est de la forme $a^2+b^2$, avec $a,b \in \Z$.
\end{corollaire}

\smallskip
      
  \n {\it Exemples } (la table des $\t_{t,n}'$ pour $n\leq 6$ et $0\leq t
\leq 2^{n-1}$ est sur le site \cite{siteweb}) :
  
  $\t_{0,1}' = \D;$
 
 $\t_{0,2}' = \D; \ \t_{1,2}' = \D^5;$
 
 $\t_{0,3}' = \D; \ \t_{1,3}' = \D^5 + \D^{13}+\D^{21} ; \ \t_{2,3}' =\D^{17} \ ; \ \t_{3,3}' = \D^{13}+\D^{21}$.
 
 \medskip

 \medskip
 \n {\it Action des opérateurs de Hecke sur les $\t_{t,n}'$}.
 
 \smallskip
 
 Si $p \equiv 3 \ {\rm ou} \ 7 \pmod 8$, on a $T_p|\t_{t,n}' = 0.$ 
 
 Si $p \equiv 1 \ {\rm ou} \ 5 \pmod 8$, on écrit $p$  sous la forme
 $p = a^2 + 4b^2$, avec  $a,b \in \Z$; on pose $t(p)' \equiv b/a  \pmod {2^n}$
 et $t^*(p)' = - t_1(p)'$.  On a:
 
 \smallskip
 
  $T_p|\t_{t,n}' = \t_{t\bullet' t(p)',n}' + \t_{t\bullet' t^*(p)',n}'$ 
 
 \smallskip
 
 \n où l'on a noté $x\bullet' y$ la loi de composition sur $\Z/2^n\Z$ définie par la formule $x\bullet' y = (x+y)/(1-4xy)$.

\n On a en particulier $\;\t'_{2^{n-1}-t(p)',n} = T_p|\D^{1+2^{2n}}$.




\end{document}